\def\Z{\mathbb Z}
  
  \def\F{\mathbb F}
  
  \def\N{\mathbb N}
 \def\R{\mathbb R}
  \def\a{\alpha}
  \def\b{\beta}

  \def\e{\epsilon}
  \def\la{\lambda}
  \def\d{\delta}

  \def\O{\Omega}
  
  \def\no{\noindent}
  \def\pf{{\it Proof. }$\;\;$}
  \def\hal{\unskip\nobreak\hfil\penalty50\hskip10pt\hbox{}\nobreak
  \hfill\vrule height 5pt width 6pt depth 1pt\par\vskip 2mm}

    \documentclass[11pt,a4paper]{article}
            \usepackage{amsfonts}

  \newtheorem{thm}{Theorem}[section]
  \newtheorem{prop}[thm]{Proposition}
  \newtheorem{lem}[thm]{Lemma}
  \newtheorem{cor}[thm]{Corollary}

\begin{document}

 \title{Rapid expansion in finite simple groups}
  \author{Martin W. Liebeck,  Gili Schul,  Aner Shalev }

  \maketitle

  \footnotetext{
The first and third authors acknowledge the support of EPSRC
Mathematics Platform grant EP/I019111/1.
The second and third authors acknowledge the support of an
ERC advanced grant 247034 and of an Israel Science Foundation grant 1117/13.
The third author acknowledges the support of the Vinik Chair
of Mathematics which he holds.}
  \footnotetext{2010 {\it Mathematics Subject Classification:}
20D06, 20F69}

\vspace{4mm}
  \begin{abstract}
We show that small normal subsets $A$ of finite simple groups expand
very rapidly --  namely, $|A^2| \ge |A|^{2-\e}$, where $\e >0$ is
arbitrarily small.
 \end{abstract}


  \section{Introduction}

In recent years there has been intense interest in the expansion of powers
of subsets of finite simple groups. For example, the remarkable product
theorem of \cite{BGT, PS} states that if $G$ is a simple group of Lie type,
and $A$ is any subset generating $G$, then either $A^3 = G$ or
$|A^3| \ge |A|^{1+\e}$, where $\e>0$ depends only on the rank of $G$
(see also \cite{helf1, helf2} for the groundbreaking results on
$L_2(p)$ and $L_3(p)$).

The case where the subset $A$ is a conjugacy class and $G$ is an arbitrary
finite simple group was considered in \cite{Sh} before the product
theorem was established.
Theorem 2.7 of \cite{Sh} shows that for any $\d>0$ there is $\e>0$ depending
on $\d$ such that $|A^3| \ge |A|^{1+\e}$ for any class $A$ of size at most
$|G|^{1-\d}$; here $G$ is any finite simple group, and $\e$ does not depend
on its rank or degree.

While the above mentioned results establish 3-step expansion,
results on $2$-step expansion were also obtained.
In \cite[10.4]{Sh} it is shown that if $A$ is a conjugacy class
of a finite simple group $G$ of Lie type, then $|A^2| \ge |A|^{1+\e}$
where $\e > 0$ now depends on the rank of $G$.
This was recently extended in \cite[1.5]{Gill} as follows: there are
absolute constants $b \in \N$ and $\e>0$ such that for any normal subset
$A$ of a finite simple group $G$, either $A^b = G$ or
$|A^2|\ge |A|^{1+\e}$.

In this paper we obtain a stronger expansion result for normal subsets -- that is, subsets which are
closed under conjugation -- as follows.

\begin{thm}\label{main} Given any $\e > 0$, there exists $b \in \N$ such that
for any normal subset $A$ of any finite simple group $G$, either $A^b = G$
or $|A^2| \ge |A|^{2-\e}$.
\end{thm}

Obviously $|A^2| \le |A|^2$, so the result says that small normal subsets of
simple groups expand almost as fast as possible. Indeed, Theorem \ref{main}
follows from

\begin{thm} \label{small}
Given any $\e > 0$, there exists $\d>0$ such that if $A$ is a normal subset
of a finite simple group $G$ satisfying $|A| \le |G|^\d$, then
$|A^2|\ge |A|^{2-\e}$.
\end{thm}

Note that some upper bound on the size of $A$ is needed in order for the
conclusion to be true.

Theorem \ref{small} holds vacuously for simple groups of bounded order
or of bounded rank,
since for these we may choose $\d$ so small that $|A| > |G|^{\d}$ for all
nontrivial classes (see Lemma \ref{vac} below); in particular, it holds for
the sporadic groups and the exceptional groups of Lie type.
It therefore remains to prove it for classical groups and alternating
groups.

We will deduce Theorem \ref{small} from the following more general result.

\begin{thm}
\label{main1}
Given any $\e > 0$, there exists $\d>0$ such that if $A_1, A_2$ are normal
subsets of a finite simple group $G$ satisfying $|A_i| \le |G|^\d$
for $i=1,2$, then $|A_1 A_2|\ge (|A_1|\,|A_2|)^{1-\e}$.
\end{thm}

Our proof of Theorem \ref{main1} is based on
results from \cite{lish99, lishdiam,  lishchar}, together with some
new results on the size of the conjugacy classes in classical groups
and in symmetric groups; see e.g. Proposition \ref{bds} below.
In fact, under the assumptions of Theorem \ref{main1}, we establish
a stronger conclusion: there exists a single conjugacy class
$C \subseteq A_1 A_2$ such that $|C|\ge (|A_1|\,|A_2|)^{1-\e}$.
The notion of the support of elements of $G$ plays a key
role in our argument.

A similar result for $k$ subsets follows inductively from
Theorem \ref{main1}:

\begin{cor}
\label{cor1}
Given $ \e> 0$ and $k \in \N$ there exists $\d > 0$ such that
if $A_1, \ldots , A_k \subseteq G$ are normal subsets of a finite simple group
$G$ with $|A_i| \le |G|^\d$, then
$|A_1 \cdots A_k|\ge (|A_1|\cdots |A_k|)^{1-\e}$.
In particular, $|A^k| \ge |A|^{k-\e}$ for every normal subset $A$ of $G$
satisfying $|A| \le |G|^\d$, where $\d$ depends on $\e$ and $k$.
\end{cor}

We also prove a result analogous to Theorem \ref{main1} for algebraic groups
over algebraically closed fields:

\begin{thm}
\label{alg}
Given any $\e > 0$, there exists $\d>0$ such that if $A_1, A_2$ are conjugacy
classes in a simple algebraic group $G$ satisfying $\dim A_i \le \d \dim G$
for $i=1,2$, then the product $A_1A_2$ contains a conjugacy class of
dimension at least $(1-\e)(\dim A_1+\dim A_2)$.
\end{thm}

The layout of the paper is as follows. In Section 2 we reduce
Theorem \ref{main1} to the case where the subsets $A_1, A_2$ are
conjugacy classes. The main result of that section is Theorem
\ref{bigclass} below, showing that (non-empty) normal subsets of classical
groups of large rank or of alternating groups of large degree
contain a relatively large conjugacy class.
Section 3 is devoted to classical groups of large rank. We study
the size of conjugacy
classes in these groups, and show that it is closely related
to the support of the elements in the class; see Propositions
\ref{bds} and \ref{eps} for more details. Section 3 concludes
the proof of Theorem \ref{main1} for classical groups.
Then, in Section 4, we prove Theorem \ref{main1} for alternating
groups, and derive some stronger results.
Finally in Section 5 we deduce Theorems \ref{main} and \ref{alg}
as well as Corollary \ref{cor1}.


\section{Reduction to conjugacy classes}

We start with some notation.
Throughout, finite simple groups $G$ are assumed to be nonabelian, and for
subsets $A_1, \ldots , A_k$ of $G$ we define
$A_1 \cdots A_k = \{ a_1 \cdots a_k : a_i \in A_i \}$.
A subset $A \subseteq G$ is said to be {\it normal} if it is closed under
conjugation, namely it is a union of conjugacy classes.
We define the {\it rank} of a finite simple group to be its untwisted
Lie rank if it is a group of Lie type, and to be its degree if it is
an alternating group.

For a finite group $G$, and a positive integer $i$, define $c_i(G)$ to be
the number of conjugacy classes of $G$ of size $i$.
For $s \in \R$, the function
\[
\eta^G(s) = \sum_{i \in \N} c_i(G)i^{-s} = \sum_C |C|^{-s}
\]
(where the second sum is over conjugacy classes $C$), was defined in
\cite{lishchar} and studied for simple groups $G$.

We start with two results which are of independent interest.

\begin{prop}\label{size}
For any $\e > 0$ there exists $N$ such
that if $G$ is a finite simple group of rank at least $N$, then for all
$m\in \N$,
$G$ has at most $m^{\epsilon}$  conjugacy classes of size at most $m$.
\end{prop}

\pf The alternating case is covered in the proof of  \cite[2.3]{lishdiam}.
So now assume that $G$ is a classical group, and let $\e>0$.
Theorem 1.10(ii) of \cite{lishchar} shows that $\eta^G(\e/2) \rightarrow 1$
as ${\rm rank}(G) \rightarrow \infty$.
Hence there exists $N$ such that for $G$ of rank at least $N$ we have
$\sum_{i \ge 1} c_i(G)i^{-\e/2} \le 1 + \e/2$, and it follows that
$\sum_{i=1}^m c_i(G) \le (1+\e/2)m^{\e/2}$. For $m \ge 3$ we have
$1 + \e/2 \le e^{\e/2} \le m^{\e/2}$, which implies
$\sum_{i=1}^m c_i(G) \le m^\e$ as required. Finally, the last inequality holds
trivially for $m=1,2$ (since $c_2(G)=0$).
The conclusion follows.
\hal
\medskip

The next result shows that normal subsets of finite simple groups
of large rank contain a relatively large conjugacy class.

\begin{thm}\label{bigclass}
For any $\e > 0$ there exists $N$
such that if $G$ is a finite simple group of rank at least $N$, and $A$
is a non-empty normal subset of $G$, then $A$ contains a conjugacy class $C$
such that $|C| \ge |A|^{1-\e}$.
\end{thm}

\pf Let $\e>0$, and let $N$ be as in the conclusion of Proposition \ref{size}.
Let $G$ be a finite simple group of rank at least $N$, and $A$ a normal
subset of $G$. Denote by $m$ the maximal size of a conjugacy class contained
in $A$. Then $A$ is a union of at most $m^\e$ classes, each of size at most $m$, and hence $|A| \le m^{1+\e}$. This implies that $m\ge |A|^{1-\e}$, and the result follows.
\hal
\medskip

\vspace{2mm}
Note that Theorem \ref{bigclass} improves \cite[2.4]{lishdiam} in the case
where the rank is unbounded.


\vspace{4mm}
We now reduce Theorem \ref{main1} to the case where the normal subsets
in the theorem are single conjugacy classes.
First we need

\begin{lem}\label{vac}
For any $N \in \N$, there exists $\d>0$ such that if a finite simple group
$G$ has a nontrivial conjugacy class of size at most $|G|^\d$,
then ${\rm rank}(G) \ge N$.
\end{lem}

\pf The case of alternating groups is trivial, since the order of the group
is then bounded in terms of the rank. Now suppose $G = G(q)$ is of Lie type
over $\F_q$ of rank $r$. Since $|x^G| = |G:C_G(x)|$, the size of a nontrivial
conjugacy class in $G$ is at least the minimal index of a proper subgroup,
which is at least $cq^r$,
where $r = {\rm rank}(G)$ and $c>0$ is a constant, as can be seen from
\cite[Tables 5.2A, 5.3A]{KL}. The result follows since $|G| < q^{4r^2}$.
\hal
\medskip

\begin{lem}\label{red}
It suffices to prove Theorem $\ref{main1}$ in the case where $A_1, A_2$
are single conjugacy classes.
\end{lem}

\pf Assume the conclusion of Theorem $\ref{main1}$ holds in the case of
conjugacy classes.
Namely, given $\e>0$, there exists $\d_1>0$ such that if $C_1,C_2$ are
conjugacy classes of finite simple group $G$ of size at most $|G|^{\d_1}$,
then $|C_1 C_2|\ge (|C_1|\,|C_2|)^{1-\e/2}$.

Applying Theorem \ref{bigclass}, choose $N$ such that whenever $A$ is a normal
subset of a simple group $G$ of rank at least $N$, then $A$ contains
a conjugacy class $C$ such that $|C| \ge |A|^{1-\e/2}$.

By Lemma \ref{vac}, there exists $\d_2>0$ such that if a finite simple group
$G$ has a nontrivial conjugacy class of size at most $|G|^{\d_2}$, then
${\rm rank}(G) \ge N$. Define $\d = {\rm min}(\d_1,\d_2)$.

Let $G$ be a finite simple group, and let $A_1, A_2$ be normal subsets of
$G$ satisfying $|A_i| \le |G|^\d$ for $i=1,2$.
Let $C_i$ be a largest conjugacy class in $A_i$, so that
$|C_i| \ge |A_i|^{1-\e/2}$. Then
\[
|A_1A_2| \ge |C_1C_2| \ge (|C_1|\,|C_2|)^{1-\e/2}  \ge
(|A_1|\,|A_2|)^{(1-\e/2)^2}  \ge (|A_1|\,|A_2|)^{1-\e},
\]
as required.
\hal
\medskip

\bigskip


\section{Classical groups}

In this section we relate the size of a conjugacy class in a
classical group to the support of the elements in the class;
we then use our results to prove Theorem \ref{main1} for
classical groups. By Lemma \ref{vac}, we need only prove the result
for classical groups of large dimension.

Let $G$ be one of the  classical groups
$L_n^\pm(q)$, $PSp_n(q)$ or $P\O_n^\pm(q)$, and let $V = V_n(q^u)$ be the natural module for $G$ with $n$ large, where $u=2$ if $G$ is unitary and $u=1$ otherwise.
Let $\bar{\F}$ be the algebraic closure of $\F_{q}$, and let $\bar{V} = V \otimes \bar{\F}$. Let $x \in G$, and let
$\hat{x}$ be a preimage of $x$ in $GL(V)$. Define
\[
\nu(x)  = \nu_{V,\bar{\F}}(x) =
\;\hbox{min}\{\dim [\bar{V},\la \hat{x}] : \la \in \bar{\F}^* \}.
\]
We shall refer to $\nu(x)$ as the {\it support} of $x$.

The following proposition, which is an extension of \cite[3.4]{lish99}, shows that $\nu(x)$  is closely related to the size of the conjugacy class of $x$. Define
\[
a(G) = \left\{ \begin{array}{l} 1, \hbox{ if }G = L_n^\pm (q) \\
                                             \frac{1}{2}, \hbox{ otherwise}
\end{array} \right.
\]

\begin{prop}\label{bds}
Suppose that $\nu(x) = s < \frac{n}{2}$, and let $a = a(G)$. There are absolute constants $c,c'>0$ such that
\[
cq^{2as(n-s-1)} \le |x^G| \le  cq^{as(2n-s+1)}.
\]
\end{prop}

\pf In the case where $x$ has prime order this is  \cite[3.4]{lish99}, but the general case requires quite a bit more argument.

Write $\hat x = tu$, where $t$ is the semisimple part and $u$ the unipotent part.

First suppose that $G = L_n(q)$.
Since $\nu (t) \le \nu(x) = s < \frac{n}{2}$, the semisimple part $t$ has an eigenvalue $\la \in \bar \F$ of multiplicity $n-s > \frac{n}{2}$. Then $\la$ must lie in $\F_q^*$. Denote by $V_\la$ the $\la$-eigenspace of $t$, and let $u$ act on $V_\la$ as $\sum_i J_i^{n_i}$, where $J_i$ is a Jordan block of size $i$. Writing $k = n-\sum in_i$, we have
\[
\hat x = \la \sum J_i^{n_i} \oplus K = x_1 \oplus K,
\]
where $x_1 =  \la \sum J_i^{n_i} \in GL_{n-k}(q)$, $K \in GL_k(q)$, and $n = k+\sum in_i = s+\sum n_i$. If we write
\begin{equation}\label{fdef}
f = \sum_i in_i^2 + 2\sum_{i<j}in_in_j,
\end{equation}
then $|C_{GL_{n-k}(q)}(x_1)| \sim q^f$ (see \cite[3.1]{LSbk}), and hence
\begin{equation}\label{ineq}
cq^{n^2-f-k^2} < |x^G| < c'q^{n^2-f-k},
\end{equation}
where $c,c'>0$ are constants. Now $\nu(x_1) = n-k-\sum n_i =s-k$. So the inequalities labelled (1) and (2)  in the proof of \cite[3.4(i)]{lish99} show that
\[
(n-s)^2+s-k \le f \le (n-s)^2+(s-k)^2.
\]
Putting this into (\ref{ineq}) gives the conclusion of the lemma for the case $G = L_n(q)$.

Next consider the unitary group $G = U_n(q)$. Again the semisimple part $t$ has an eigenspace $V_\la$ of dimension greater than $\frac{n}{2}$. Write $(\,,\,)$ for the unitary form on $V$ preserved by $G$, and $\a \rightarrow \bar \a$ for the involutory automorphism of the field $\F_{q^2}$. There is a nonsingular vector $v \in V_\la$, so $0 \ne (v,v) = (v\hat x, v\hat x) = \la \bar{\la} (v,v)$, and hence $\la \bar \la = 1$. Also $V_\la$ is a non-degenerate subspace, since its radical must be contained in the radical of the whole space $V$. Hence, letting $u$ acts on $V_\la$ as $\sum_i J_i^{n_i}$, we have as above
\[
\hat x = \la \sum J_i^{n_i} \perp K,
\]
where $K \in GU_k(q)$ and $n = k+\sum in_i = s+\sum n_i$. Now we argue exactly as in the previous paragraph.

Next let $G = PSp_n(q)$. Here $t$ has an eigenspace $V_\la$ of dimension greater than $\frac{n}{2}$, and $\la$ must be $\pm 1$ and $V_\la$ non-degenerate. So $\hat x =  \la \sum J_i^{n_i} \perp K = x_1\perp K$ with $K \in Sp_k(q)$ and $n = k+\sum in_i = s+\sum n_i$.

Suppose $q$ is odd. Then by \cite[3.1]{LSbk} we have $|C_{Sp_{n-k}(q)}(x_1)| \sim q^{g/2}$, where
\begin{equation}\label{gdef}
g =  \sum_i in_i^2 + 2\sum_{i<j}in_in_j + \sum_{i\;odd}n_i.
\end{equation}
As $|Sp_n(q)|\sim q^{\frac{1}{2}(n^2+n)}$, it follows that
\begin{equation}\label{ineq1}
cq^{\frac{1}{2}(n^2+n-g-k^2-k)} < |x^G| < c'q^{\frac{1}{2}(n^2+n-g-k)}.
\end{equation}
The inequalities labelled (1) and (2)  in the proof of \cite[3.4(ii)]{lish99} show that
\[
(n-s)^2+n-s \le g \le (n-s)^2+(s-k)^2 + n-k.
\]
Putting this into (\ref{ineq1}) gives the conclusion of the lemma for $G = PSp_n(q)$ with $q$ odd.

Now suppose $q$ is even. This is slightly more complicated, as in general there can be many unipotent classes in a symplectic group having the same Jordan form. The general form of a unipotent element, and its centralizer, is given by \cite[6.2, 7.3]{LSbk}, from which it can be seen that $|C_{Sp_{n-k}(q)}(x_1)| \sim q^{\frac{1}{2}g'}$, where
\[
g\le g' \le g+2\sum_{i,n_i\;even}n_i
\]
and $g$ is as above. Then $g'\ge g \ge (n-s)^2+n-s$, and the lower bound for $|x^G|$ follows as before. As for the upper bound, observe that
\[
\begin{array}{ll}
(s-k)^2+n-k & = \left(\sum (i-1)n_i\right)^2 + \sum in_i \\
              & \ge \sum(i-1)n_i^2 +2\sum_{i<j} (i-1)n_in_j + \sum_{i\;odd}n_i + 2\sum_{i\;even}n_i \\
               & = g + 2\sum_{i\;even}n_i - (n-s)^2 \\
               & \ge g' - (n-s)^2.
\end{array}
\]
Hence $g' \le (n-s)^2 +(s-k)^2+n-k$, and the upper bound for $|x^G|$ follows as before.
This completes the proof for the symplectic groups.

The argument for orthogonal groups is very similar: again
we have $\hat x  =  \la \sum J_i^{n_i} \perp K = x_1\perp K \in O_{n-k}(q) \times O_k(q)$, where $\la = \pm 1$, $k<\frac{n}{2}$ and $n = k+\sum in_i = s+\sum n_i$. If we define
\[
h =  \sum_i in_i^2 + 2\sum_{i<j}in_in_j - \sum_{i\;odd}n_i,
\]
then for $q$ odd we have $|C_{O_{n-k}(q)}(x_1)| \sim q^{h/2}$, and for $q$ even we have
$|C_{O_{n-k}(q)}(x_1)| \sim q^{h'/2}$, where $h-2\sum_{i,n_i\;even}n_i \le h'\le h$ (see \cite[3.1, 6.2, 7.3]{LSbk}). Arguing as for the symplectic case, we see that
\[
(n-s)^2 -(n-k)-(s-k) \le h' \le (n-k)^2+(s-k)^2+2(s-k)-n,
\]
and the conclusion follows. \hal

\begin{lem}\label{bigs}
Let $x \in G$ with $\nu(x) = s$, and suppose that $|x^G| \le |G|^{\frac{1}{4}}$. Then $s < \frac{n}{2}-1$.
\end{lem}

\pf First suppose $G = L_n(q)$, and write $\hat x = tu$ as in the previous proof. Recall our assumption that $n$ is large. The centralizer of $t$ in $GL_n(q)$ is of the form $C = \prod GL_{n_i}(q^{a_i})$, where $\sum n_ia_i = n$. Since this must have order greater than $|G|^{\frac{3}{4}}$, it follows that the largest factor of $C$ is $GL_r(q)$, where $r>\frac{n}{2}$ and $q^{r^2+(n-r)^2} > |G|^{\frac{3}{4}}$. Hence in fact $r>\a n$, where $\a = 0.85$. Let $V_r$ be the $r$-dimensional eigenspace for $t$, and let $u$ act on $V_r$ as $\sum J_i^{n_i}$. So as in the previous proof we have
\[
\hat x = \la \sum J_i^{n_i} \oplus K = x_1 \oplus K \in GL_r(q) \times GL_{n-r}(q).
\]
Let $s_1 = \nu(x_1)$, so that $s \le s_1+n-r$. Define $f$ as in (\ref{fdef}) in the previous proof.

Suppose $s_1>\frac{r}{2}$. Then the inequality (3) in the proof of \cite[3.4(i)]{lish99} shows that $f\le r(r-s_1)$. Therefore $|x^G| \ge |x_1^{GL_r(q)}| \ge cq^{rs_1}$ (where $c$ is a positive constant). Since by hypothesis $|x^G|\le |G|^{\frac{1}{4}}$, it follows that $rs_1 \le \frac{n^2}{4}$. Then
\begin{equation}\label{lnc}
s \le s_1+n-r \le \frac{n^2}{4r}+n-r,
\end{equation}
which is less than $\frac{n}{2}-1$ since $r > \a n$.

Now suppose $s_1 \le \frac{r}{2}$. Then the inequality (2) in the proof of \cite[3.4(i)]{lish99} shows that
$f\le (r-s_1)^2+s_1^2$, and so $|x^G| \ge |x_1^{GL_r(q)}| \ge cq^{2s_1(r-s_1)}$.
Thus $2s_1(r-s_1) \le \frac{n^2}{4}$. Writing $\b = \frac{s_1}{r}$ (so $0<\b<\frac{1}{2}$), this gives $2\b(1-\b)r^2 \le \frac{n^2}{4}$, and hence
\begin{equation}\label{ab}
2\b(1-\b) \le \frac{1}{8\a^2}.
\end{equation}
Also $s \le s_1+n-r \le n-(1-\b)r \le n(1-\a(1-\b))$. Now check that for $\b$ satisfying (\ref{ab}), we have $\a(1-\b)>\frac{1}{2}$, and the conclusion follows.
This completes the proof for $G = L_n(q)$.

The proof for the other classical groups runs along entirely similar lines. We shall  just give a sketch for the symplectic groups and leave the other cases to the reader. Let $G = PSp_n(q)$ with $n$ large, and write $\hat x = tu$ as above. The centralizer of $t$ in $Sp_n(q)$ is of the form $C = Sp_r(q) \times Sp_s(q) \times \prod GL^{\e_i}_{n_i}(q^{a_i})$, where $n = r+s+2\sum n_ia_i$ and the first two factors correspond to the $\pm 1$-eigenspaces. This has order greater than $|G|^{\frac{3}{4}}$, so $C$ must have a factor $Sp_r(q)$, where $r>\frac{n}{2}$ and $|Sp_r(q)\times Sp_{n-r}(q)| \ge |G|^{\frac{3}{4}}$. As above it follows that for large $n$ we have $r>\a n$ with $\a = 0.85$.
As usual we can write
\[
\hat x = \la \sum J_i^{n_i} \oplus K = x_1 \oplus K \in Sp_r(q) \times Sp_{n-r}(q),
\]
where $\la = \pm 1$. As in the proof of the previous lemma we have $|C_{Sp_r(q)}(x_1)| \sim q^{\frac{1}{2}g'}$, where
$g \le g'\le g+2\sum_{i,n_i\;even}$ and $g$ is as in (\ref{gdef}). If $s_1>\frac{r}{2}$ then $g'\le (r-s_1)^2+s_1(r-s_1)+r$, since
\[
\begin{array}{lll}
(r-s_1)^2+s_1(r-s_1)+r & = & \left(\sum n_i\right)^2+\left(\sum (i-1)n_i\right)\left(\sum n_i\right) + \sum in_i \\
                                  & \ge & \sum n_i^2 + 2\sum_{i<j}n_in_j + \sum (i-1)n_i^2 + \\
                                  & & 2\sum_{i<j}(i-1)n_in_j + \sum in_i \\
                                  & \ge & g + 2\sum_{i\;even}n_i \\
                                  & \ge & g'.
\end{array}
\]
Hence $|x^G| \ge |x_1^{Sp_r(q)}| \ge cq^{\frac{1}{2}(r^2+r-g')} \ge
cq^{\frac{1}{2}rs_1}$. As $|x^G|\le |G|^{\frac{1}{4}}$ it follows that $rs_1\le \frac{n^2+n}{4}$. Now the conclusion follows as in (\ref{lnc}) above. Finally, if $s_1\le \frac{r}{2}$  then we similarly deduce that $g' \le (r-s_1)^2+s_1^2+r$, which implies that $|x^G| \ge cq^{s_1(r-s_1)}$. Hence $s_1(r-s_1) \le \frac{n^2+n}{4}$ and now we argue as in the $L_n(q)$ case above. \hal

\begin{lem}\label{first}
Let $x\in G$ with $\nu(x)=s$, and let $0<\d \le \frac{1}{4}$.
There is a constant $d$ such that if $|x^G| \le |G|^\d$, then
$s \le 2\d n + \frac{d}{n}$.
\end{lem}

\pf Let $C = x^G$ and suppose $|C| \le |G|^\d$. By Lemma \ref{bigs} and
Proposition \ref{bds} we have
\[
cq^{2as(n-s-1)} \le |C| \le |G|^\d < q^{n^2\d},
\]
Writing $d' = \log_2\frac{1}{c}$, this gives
$q^{2as(n-s-1)} < q^{n^2\d +d'}$. Since $s\le \frac{n}{2}-1$,
this implies $ans < \d n^2+ d'$. \hal

The next result shows that the size of a small conjugacy class
$x^G$ of a finite simple classical group $G$ is almost determined
by the support $\nu(x)$ of $x$.

\begin{prop}\label{eps}
For any $\e_1>0$, there exists $\d>0$ such that if $x \in G$ with $\nu(x)=s$ and $|x^G|\le |G|^\d$, then
\[
q^{(2a-\e_1)ns} \le |x^G| \le q^{(2a+\e_1)ns}.
\]
\end{prop}

\pf We may assume that $\e_1< \frac{2}{3}$. Choose $\d = \frac{\e_1}{4}$. Now $s \le 3\d n$ for large $n$, by Lemma \ref{first}. Since $s<\frac{n}{2}$, we may apply Proposition \ref{bds}. We have $\e_1n \ge 3\d n+1 \ge s+1$, so Proposition \ref{bds} gives the conclusion. \hal

\vspace{4mm}
Now let $x_1,x_2 \in G$, and assume that $\nu(x_i) = s_i$ with $s_i < \frac{1}{4}n$ for $i=1,2$.
The largest eigenspace of $\hat x_i$ on $\bar V$ has dimension $n-s_i > \frac{3}{4}n$, and it follows that the corresponding eigenvalue $\la_i$ lies in $\F_{q^u}$, and also that $\la_i\bar \la_i = 1$ in the unitary case, and $\la_i = \pm 1$ in the symplectic and orthogonal cases. As in the proof of Proposition \ref{bds} we have
$\hat x = \la \sum J_i^{n_i} \perp K$, and separating the Jordan blocks of size 1, we can write
\[
\hat x_i = \la_i I_{t_i} \perp \sum_{j=1}^{r_i} J_{n_{ji}}(\la_i) \perp K_i,
\]
where $J_{n_{ji}}(\la_i)$ denotes a single Jordan block of size $n_{ji}\ge 2$ for each $j$, and $K_i$ has no eigenvalue equal to $\la_i$; moroever the subspaces on which the three summands act are non-degenerate and mutually perpendicular in the case $G \ne L_n(q)$.

Now $s_i = n-(t_i+r_i)$ and $n\ge t_i+2r_i$, hence $t _i> n-s_i-\frac{1}{2}(n-t_i)$. Since $s_i < \frac{1}{4}n$ it follows that $t_i>\frac{1}{2}n$, and
\[
\hat x_i = \la_i I_{t_i} \perp L_i,
\]
where $L_i = \sum_j J_{n_{ji}}(\la_i) \perp K_i$. Now define
\[
\hat y = \la_1\la_2 I_{t_1+t_2-n} \perp \la_2 L_1 \perp \la_1L_2,
\]
and let $y$ be the image of $\hat y$ in $G$. Write $y = x_1 * x_2$ (defined only up to conjugacy).

 \begin{lem}\label{yprop}
Let $y = x_1*x_2$ as above. Then $y \in x_1^Gx_2^G$, and $\nu(y) = \nu(x_1)+\nu(x_2)$.
\end{lem}

\pf There are conjugates of $\hat x_1,\hat x_2$ of the form $ \la_1 I_{t_1+t_2-n} \perp L_1 \perp \la_1 I_{n-t_2}$ and
$ \la_2 I_{t_1+t_2-n} \perp \la_2 I_{n-t_1} \perp L_2$ respectively, and their product is equal to $\hat y$.
Also, from the definition of $\hat y$ we have
\[
\nu (y) = n - (t_1+t_2-n) - r_1-r_2 = 2n - (t_1+r_1) - (t_2+r_2) = s_1+s_2,
\]
as required. \hal

\begin{lem}\label{nextone}
Given $\e>0$, there exists $\d>0$ such that the following holds. If $x_1^G, x_2^G$ are classes in $G$ with $|x_i^G| \le |G|^\d$ for $i=1,2$, and $y = x_1*x_2$, then $|y^G| \ge (|x_1^G|\,|x_2^G|)^{1-\e}$.
\end{lem}

\pf By Lemma \ref{first} and Proposition \ref{eps}, there exists $\d>0$ such that if $x^G$ is a class such that $|x^G| \le |G|^\d$, then $\nu(x) < \frac{n}{4}$ and also the conclusion of Proposition \ref{eps} holds with $\e_1 = \frac{\e}{2}$.

Now let $x_1^G, x_2^G$ be classes in $G$ with $|x_i^G| \le |G|^\d$. Then $s_i:=\nu(x_i) < \frac{n}{4}$ for $i=1,2$, so we can define $y = x_1*x_2$. Moreover $\nu(y) = s_1+s_2$ by Lemma \ref{yprop}, so Proposition \ref{eps} gives
\[
|y^G| \ge q^{(2a-\e_1) n(s_1+s_2)},\;\;|x_1^G|\,|x_2^G|
\le q^{(2a+\e_1) n(s_1+s_2)}.
\]
The conclusion follows, since $\frac{2a-\e_1}{2a+\e_1} \ge 1-\e$. \hal

\vspace{4mm}
Theorem \ref{main1} for classical groups now follows from Lemmas \ref{yprop}
and \ref{nextone}, together with Lemma \ref{red}.

\section{Alternating groups}

We now prove Theorem \ref{main1} for alternating groups.
Recall that it suffices to prove it for conjugacy classes.
We start with symmetric groups $S_n$.

Let $\pi \in S_n$ and let $s$ be the support of $\pi$, namely
\[
s = |\{1 \le j \le n : \pi(j) \ne j \}|.
\]
As in the previous section, we shall relate the size of
the conjugacy class of an element $\pi \in S_n$ to the support
$s$ of $\pi$. However, in this case it is not true that the
support almost determines the class size as in Proposition
\ref{eps} for classical groups.

For $i \ge 2$ let $c_i$ denote the number of cycles of length
$i$ in $\pi$. Then $s = \sum_{i \ge 2} ic_i$ and $\pi$ has
$n-s$ fixed points.
Let $C$ be the conjugacy class of $\pi$ in $S_n$. It is well known
that
\[
|C| = \frac{n!}{(n-s)! \prod_{i\ge 2} i^{c_i} \prod_i c_i!}.
\]

\begin{lem}\label{lem:Sym size C}
With the above notation we have
$\frac{n!}{(n-s)!2^{s/2}(s/2)!} \le |C| \le \frac{n!}{(n-s)!}$.
\end{lem}

\pf The upper bound is trivial.
For the lower bound, note that $\sum_{i \ge 2} c_i \le
\frac{1}{2} \sum_{i \ge 2} ic_i = \frac{s}{2}$. This implies
$\prod_{i\ge 2} c_i! \le \frac{s}{2}!$.
Also $\prod_{i \ge 2} i^{c_i} \le \prod_{i \ge 2} 2^{ic_i/2} = 2^{s/2}$. The lower bound follows
from these inequalities. \hal

\vspace{4mm}
Note that the lower bound above is best possible, as shown by
the case $\pi = (12)(34)\ldots (s-1 s)$ ($s$ even). The upper
bound is almost best possible (take $\pi = (1 2 \ldots s)$).

We fix some notation for the rest of this section.
Let $\pi_{1},\pi_{2}\in S_n$ be permutations of supports $s_{1}, s_{2}$
respectively. For $i=1,2$ let $C_i$ denote the conjugacy class of $\pi_i$
in $S_n$.

Suppose $s_1 + s_2 \le n$. Then there exists a permutation, which we
denote by $\pi_{2}^{\prime}$, that has the same cycle structure as
$\pi_{2}$, such that the points moved by $\pi_{2}^{\prime}$ are
fixed points of $\pi_1$.
Define the conjugacy class
\begin{equation}\label{cstar}
C_1 * C_2 = (\pi_{1}\pi_{2}^{\prime})^{S_n} \subseteq C_{1}C_{2}.
\end{equation}
Note that the elements of $C_1 * C_2$ have support $s_1+s_2$.
We shall prove that $|C_1 C_2|$ is large by providing lower bounds
on the size of $|C_1 * C_2|$.

We start by showing that the conclusion of Theorem \ref{main1} holds
for conjugacy classes $C_1, C_2$ of bounded support.

\begin{lem}\label{cor:Sym constant s}
Let $s_1,s_2$ be positive integers, and let $\e>0$. There exists an integer $N=N(\epsilon,s_{1},s_{2})$
such that if $n \ge N$ and $C_1,C_2$ are classes in $S_n$ of support $s_1,s_2$ respectively, then
\[
\left|C_{1}C_{2}\right|\ge
\left(\left|C_{1}\right|\left|C_{2}\right|\right)^{1-\e}.
\]
\end{lem}

\pf We shall choose $N \ge s_1+s_2$ so the conjugacy class $C_1 * C_2$
may be constructed. Applying Lemma \ref{lem:Sym size C} for this class
(whose support is $s_1+s_2$) we obtain
\[
|C_{1}C_{2}| \ge |C_1 * C_2| \ge \frac{n!}{(n-s_{1}-s_{2})!2^{s_{1}+s_{2}}(s_{1}+s_{2})!} : = f(n).
\]
By the same lemma we also have
\begin{equation}\label{no2}
\left(\left|C_{1}\right|\left|C_{2}\right|\right)^{1-\epsilon}\le\left(\frac{n!}{\left(n-s_{1}\right)!}\cdot\frac{n!}{\left(n-s_{2}\right)!}\right)^{1-\epsilon} : = g(n)^{1-\e}.
\end{equation}
Since $f(n)$ and $g(n)$ are polynomials in $n$ of degree $s_1+s_2$, there exists $N = N(\e,s_1,s_2)$
such that $f(n) \ge g(n)^{1-\e}$ for $n \ge N$, and the conclusion follows. \hal

\vspace{4mm}
In Proposition \ref{prop:sym big n} below we will derive a similar
conclusion assuming only $s_1 , s_2 \le \frac{n}{3}$.
We need some preparations.

\begin{lem}
\label{lem:Sym term1} Suppose $s_1, s_2 \le n/3$, and define
$f_1(n,s_1,s_2) = \frac{n!}{\left(n-s_{1}-s_{2}\right)!}$ and
$f_2(n,s_1,s_2) = \frac{n!^{2}}{\left(n-s_{1}\right)!\left(n-s_{2}\right)!}$.
Then
\[
\frac{f_1(n,s_1,s_2)}{f_2(n,s_1,s_2)} \ge3^{-s_{1}}.
\]
\end{lem}

\pf We have
\[
\begin{array}{lcl}
\frac{f_1(n,s_1,s_2)}{f_2(n,s_1,s_2)} &=& \frac{\left(n-s_{1}\right)!\left(n-s_{2}\right)!}{n!\left(n-s_{1}-s_{2}\right)!}\\
&=& \frac{\left(n-s_{2}\right)\left(n-s_{2}-1\right)\cdots\left(n-s_{2}-s_{1}+1\right)}{n\left(n-1\right)\cdots\left(n-s_{1}+1\right)}\\
&\ge& \frac{\left(n-s_{1}-s_{2}\right)^{s_{1}}}{n^{s_{1}}}=\left(1-\frac{s_{1}+s_{2}}{n}\right)^{s_{1}}\ge\left(1-\frac{2}{3}\right)^{s_{1}}=\left(\frac{1}{3}\right)^{s_{1}}.
\end{array}
\]
\hal
\medskip

As above let $\pi_1,\pi_2 \in S_n$ have supports $s_1,s_2$ respectively. For $i\ge 2$ let
$c_i$ be the number of cycles of length $i$ in $\pi_1$, and $d_i$ the number
of cycles of length $i$ in $\pi_2$.

\begin{lem} \label{lem:Sym term2} We have
\[
\frac{\prod_{i\ge 2} c_{i}! \prod_{i\ge 2} d_{i}!}
{\prod_{i\ge 2}\left(c_{i}+d_{i}\right)!} \ge 2^{-(s_1+s_2)/2}.
\]
\end{lem}

\pf Observe that
\[
\frac{c_{i}!d_{i}!}{\left(c_{i}+d_{i}\right)!}=
{1 \over {{c_i+d_i}\choose{c_i}}}\ge\frac{1}{2^{c_{i}+d_{i}}}.
\]
Combining this with the inequalities $\sum c_{i} \le \sum \frac{i}{2}c_{i} =
\frac{s_{1}}{2}$ and $\sum d_i \le \frac{s_2}{2}$ we obtain
\[
\frac{\prod c_{i}! \prod d_{i}!}{\prod\left(c_{i}+d_{i}\right)!}\ge
2^{-\sum c_{i}-\sum d_{i}}\ge2^{-\frac{s_{1}+s_{2}}{2}}.
\]
\hal

\vspace{4mm}
Recall that $C_i = \pi_i^{S_n}$ for $i=1,2$.

\begin{lem}
\label{lem:Sym term3} Suppose $s_1 \le \frac{n}{2}$. Then
$|C_1| \ge s_{1}^{s_1/2}$.
\end{lem}

\pf Applying Lemma \ref{lem:Sym size C}, we obtain
\[
\begin{array}{lcl}
|C_1| = \frac{n!}{\left(n-s_{1}\right)!\cdot\prod i^{c_{i}}\cdot\prod c_{i}!}
&\ge&
\frac{n!}{\left(n-s_{1}\right)!\cdot2^{s_1/2}\left(\frac{s_{1}}{2}\right)!}\\
&\ge& \frac{\left(n-s_{1}\right)^{s_{1}}}{2^{s_1/2}\left(\frac{s_{1}}{2}\right)^{s_1/2}}\ge\frac{s_{1}^{s_{1}}}{2^{s_1/2}\left(\frac{s_{1}}{2}\right)^{s_1/2}}=s_1^{s_1/2}.
\end{array}
\]
\hal

\begin{prop}
\label{prop:sym big n} For any $\epsilon>0$ there exists
$N=N\left(\epsilon\right)$ such that if $n \ge N$ and $s_1, s_2 \le \frac{n}{3}$,
then $\left|C_{1}C_{2}\right|\ge\left(\left|C_{1}\right|\left|C_{2}\right|\right)^{1-\epsilon}$.
\end{prop}

\pf We have
\[
\begin{array}{l}
|C_1| = \frac{n!}{(n-s_1)! \prod i^{c_i} \prod c_i!},\; |C_2| = \frac{n!}{(n-s_2)! \prod i^{d_i} \prod d_i!}, \hbox{ and }\\
|C_1*C_2| = \frac{n!}{(n-s_1-s_2)! \prod  i^{c_i+d_i} \prod (c_i+d_i)!}.
\end{array}
\]
Hence, taking $f_1,f_2$ as in Lemma \ref{lem:Sym term1},
\[
\frac{|C_1*C_2|}{|C_1|\,|C_2|} = \frac{f_1(n,s_1,s_2)}{f_2(n,s_1,s_2)} \frac{\prod \left(c_{i}!\right)\prod \left(d_{i}!\right)}{\prod \left(c_{i}+d_{i}\right)!}.
\]
It follows using Lemmas  \ref{lem:Sym term1}, \ref{lem:Sym term2} and \ref{lem:Sym term3} that
\[
\begin{array}{ll}
\frac{|C_1C_2|}{(|C_1|\,|C_2|)^{1-\e}} & =  \frac{|C_1C_2|}{|C_1|\,|C_2|} |C_1|^\e |C_2|^\e \\
                                             & \ge  \frac{|C_1*C_2|}{|C_1|\,|C_2|} |C_1|^\e |C_2|^\e \\
                                             & \ge 3^{-s_{1}}2^{-(s_{1}+s_{2})/2}s_{1}^{s_1\e/2}  s_{2}^{s_2\e/2}.
\end{array}
\]

Let $S_{0}=S_{0}\left(\epsilon\right)$ be such that
$3^{-s_{1}}2^{-(s_{1}+s_{2})/2} s_{1}^{s_1\e/2}  s_{2}^{s_2\e/2} \ge 1$ provided $s_1 \ge S_{0}$ or $s_2 \ge S_0$.
In that case we deduce that
$\left|C_{1}C_{2}\right|\ge\left(\left|C_{1}\right|\left|C_{2}\right|\right)^{1-\epsilon}$.
Otherwise we have $s_1, s_2 \le S_{0}$;
let $N_0 = N_{0}\left(\epsilon,S_{0},S_{0}\right)=N\left(\epsilon\right)$
be such that for $n \ge N_{0}$, $\left|C_{1}C_{2}\right|\ge\left(\left|C_{1}\right|\left|C_{2}\right|\right)^{1-\epsilon}$
(such $N_{0}$ exists by Lemma \ref{cor:Sym constant s}). \hal

\vspace{4mm}
We will need the following well known Stirling approximation which holds for all $n$ (see \cite[\S 2.9]{feller}):
\begin{equation}\label{stirling}
\sqrt{2\pi n}\left(\frac{n}{e}\right)^{n}\le n!\le2\sqrt{2\pi n}\left(\frac{n}{e}\right)^{n}.
\end{equation}

The next result provides a lower bound close to $n^{s/2}$
for the size of a conjugacy class of support $s$.

\begin{lem}
\label{lem:Sym size C stirling}If $C\subseteq S_n$
is a conjugacy class with support $s$, then
\[
\left|C\right|\ge\frac{1}{4\sqrt{\pi n}}n^{\frac{s}{2}}e^{-\frac{s}{2}}.
\]
\end{lem}

\pf Using Lemma \ref{lem:Sym size C} and (\ref{stirling}), we obtain
\[
\begin{array}{lcl}
\left|C\right|
{\ge}\frac{n!}{\left(n-s\right)!2^{\frac{s}{2}}\left(\frac{s}{2}\right)!}
&{\ge}& \frac{\sqrt{2\pi n}}{4\sqrt{2\pi\left(n-s\right)\cdot2\pi\frac{s}{2}}}\cdot\frac{n^{n}}{\left(n-s\right)^{\left(n-s\right)}\left(\frac{s}{2}\right)^{\frac{s}{2}}}e^{-\frac{s}{2}}2^{-\frac{s}{2}}\\
&=& \frac{\sqrt{n}}{4\sqrt{\pi\left(n-s\right)s}}\cdot\frac{n^{n}}{\left(n-s\right)^{\left(n-s\right)}s^{\frac{s}{2}}}e^{-\frac{s}{2}}\ge\frac{\sqrt{n}}{4\sqrt{\pi n^{2}}}\cdot\frac{n^{n}}{n^{\left(n-s\right)}n^{\frac{s}{2}}}e^{-\frac{s}{2}}\\
&=&\frac{1}{4\sqrt{\pi n}}n^{\frac{s}{2}}e^{-\frac{s}{2}},
\end{array}
\]
as required. \hal

\begin{prop}
\label{prop:Sym final}
For any $\e > 0$ there exists $N = N(\e)$ such that, if $n \ge N$,
and $C_1, C_2$ are conjugacy classes of $G = S_n$ satisfying
$|C_1|, |C_2| \le |G|^{1/8}$,
then
$|C_1 C_2| \ge (|C_1| |C_2|)^{1-\e}$.

In particular, Theorem $\ref{main1}$ holds for conjugacy classes
in $S_n$.
\end{prop}

\pf
We will first show that if $n \ge 40$ and $C \subset G=S_n$ is a class
of support $s$ satisfying $\left|C\right|\le\left|G\right|^{\frac{1}{8}}$,
then $s\le\frac{n}{3}$.

Let $C$ be such a class. By Lemma \ref{lem:Sym size C stirling} and (\ref{stirling}),
\[
\frac{1}{4\sqrt{\pi n}}n^{\frac{s}{2}}e^{-\frac{s}{2}}\le\left|C\right|\le\left|G\right|^{\frac{1}{8}}=\left(n!\right)^{\frac{1}{8}}\le
\left(2\sqrt{2\pi n}\cdot n^{n}e^{-n}\right)^{\frac{1}{8}}.
\]
Thus
\[
1\le2^{\frac{35}{16}}\left(\pi n\right)^{\frac{9}{16}}\left(\frac{n}{e}\right)^{\frac{n}{8}-\frac{s}{2}}.
\]
Suppose $s>\frac{n}{3}$. Then $\frac{n}{8}-\frac{s}{2}<\frac{n}{8}-\frac{n}{6}=-\frac{n}{24}$.
Therefore
\[
1\le 2^{\frac{35}{16}}\left(\pi n\right)^{\frac{9}{16}}\left(\frac{n}{e}\right)^{-\frac{n}{24}}
< 9 n^{\frac{9}{16}} \left(\frac{n}{e}\right)^{-\frac{n}{24}},
\]
which is a contradiction for $n\ge 40$.
Hence $s\le \frac{n}{3}$.

Now let $C_1,C_2$ be classes as in the statement of the proposition, with supports $s_1,s_2$. By the above $s_1,s_2 \le \frac{n}{3}$, so we can take $N=N\left(\epsilon\right)$  as in Proposition \ref{prop:sym big n}.
Then if $n\ge\max\left\{ N,40\right\} $ we have $\left|C_{1}C_{2}\right|\ge\left(\left|C_{1}\right|\left|C_{2}\right|\right)^{1-\epsilon}$.
This completes the proof. \hal

\vspace{4mm}
We finally turn to alternating groups, proving

\begin{lem}
Proposition $\ref{prop:Sym final}$ holds for conjugacy classes in $G=A_n$.
\end{lem}

\pf
If $\pi\in A_n$ then $\pi^{A_n}=\pi^{S_n}$
or $\left|\pi^{A_n}\right|=\frac{1}{2}\left|\pi^{S_n}\right|$.
So the proof of Proposition \ref{prop:Sym final}
can be adjusted to $A_n$,
just by taking $n \ge 45$ instead of $40$ to get $s\le\frac{n}{3}$.
For $\pi_i\in S_n$ with support $s_i \le \frac{n}{3}$, we have
$\pi_i^{A_n}=\pi_i^{S_n}$.
Also for $\sigma=\pi_1\pi_2^{\prime}$ as in (\ref{cstar}), the support
is at most
$\frac{2n}{3}$ and so $\sigma^{A_n}=\sigma^{S_n}$,
and the proof of the proposition can continue as it is (replacing $40$ with
$45$).
\hal

\vspace{4mm}
This completes the proof of Theorem \ref{main1} for alternating
groups $G=A_n$ in a somewhat stronger form: it suffices to assume
that the normal subsets have size at most $|G|^{1/8}$ and that
$n \ge N(\e)$.

\section{Final deductions}

\vspace{4mm}
\no {\bf Deduction of Theorem \ref{main}}

\vspace{2mm}
Let $\e>0$ and let $\d>0$ be as in the conclusion of Theorem \ref{small}.
Theorem 1.1 of \cite{lishdiam} states that there is an absolute constant
$c$ such that  for every nontrivial normal subset $A$ of a finite simple
group $G$, we have $A^m = G$ for any $m \ge c\frac{\log |G|}{\log |A|}$.
Define $b = \lceil \frac{c}{\d} \rceil$.

Now let $A$ be a normal subset of a finite simple group $G$.
If $|A| \ge |G|^\d$ then the previous paragraph shows that $A^b = G$.
Otherwise, Theorem \ref{small} shows that $|A^2| \ge |A|^{2-\e}$.
This completes the proof. \hal

\vspace{4mm}
\no {\bf Deduction of Corollary \ref{cor1}}

\vspace{2mm}
We argue by induction on $k \ge 2$. The case $k=2$ is Theorem \ref{main1}.
Suppose $k \ge 3$. By induction, given $\e>0$ and $2 \le m < k$,
there exists $\d(\e,m) > 0$ such that if $A_1, \ldots , A_m \subseteq G$
are normal subsets with $|A_i| \le |G|^{\d(\e,m)}$
then $|A_1 \cdots A_m| \ge (|A_1| \cdots |A_m|)^{1-\e}$.

Define $\d(\e, k) = \min \{ \d(\e/2,2)/(k-1), \d(\e/2, k-1) \}$.

Now let $\d = \d(\e, k)$ and suppose $A_1 , \ldots , A_k$ are normal subsets
of $G$ of size at most $|G|^\d$.
By induction it follows that
\[
|A_1 \cdots A_{k-1}| \ge  (|A_1| \cdots |A_{k-1}|)^{1-\e/2}.
\]
Note that
$|A_1 \cdots A_{k-1}| \le |G|^{(k-1)\d} \le |G|^{\d(\e/2,2)}$,
and so the case $k=2$ yields
\[
|A_1 \cdots A_k| \ge (|A_1 \cdots A_{k-1}| |A_k|)^{1-\e/2} \ge
((|A_1| \cdots |A_{k-1}|)^{1-\e/2}|A_k|)^{1-\e/2},
\]
which is at least $(|A_1| \cdots |A_k|)^{1-\e}$.
The result follows.
\hal

\vspace{4mm}
\no {\bf Deduction of Theorem \ref{alg}}

\vspace{2mm}
The proof is virtually the same as that of Theorem \ref{main1}.
As in Lemma \ref{vac}, since every conjugacy class in
a simple algebraic group $G$ has dimension at least $r = {\rm rank}(G)$,
we need only consider classical groups of large dimension.
So let $G = SL_n(K)$, $Sp_n(K)$ or $SO_n(K)$ where
$K$ is algebraically closed and $n$ is large, and define $a:=a(G) = 1$,
$\frac{1}{2}$ or $\frac{1}{2}$,
respectively. Let $x \in G$ and define $s = \nu(x)$ as in Section 3.
The proof of Lemma \ref{bds} gives
\[
2as(n-s-1) \le \dim x^G \le as(2n-s+1),
\]
and there are similar dimensional analogues of Lemma \ref{bigs},
Lemma \ref{first} and Proposition \ref{eps},
with the same proofs. For $x_1,x_2 \in G$ with
$s_i = \nu(x_i) < \frac{n}{4}$ we can define $y = x_1*x_2$
as before, and $\nu(y) = s_1+s_2$ as in Lemma \ref{yprop}.
Now the theorem follows as in the proof of Lemma \ref{nextone}.
\hal

\vspace{1cm}
\no M. W. Liebeck, Department of Mathematics,  Imperial College, London SW7 2AZ, UK

\no G. Schul,  Institute of Mathematics, Hebrew University, Jerusalem 91904,
Israel

\no A. Shalev,  Institute of Mathematics, Hebrew University, Jerusalem 91904,
Israel

\end{document}